\documentclass[a4paper]{article}

\usepackage[english]{babel}
\usepackage[utf8]{inputenc}
\usepackage{times}
\usepackage{pifont}
\usepackage{graphicx}
\usepackage{color}
\usepackage{amsmath}
\usepackage{amsthm}
\usepackage{amsfonts}
\usepackage{amssymb}
\usepackage{mathtools}
\usepackage{url}

\title{Distributions of $n$th Powers in Finite Fields}

\author{Aaron Doman}

\date{}

\begin{document}
\maketitle
\begin{abstract}
In this paper, we first find the distribution of $n$th power residues modulo a prime $p$ by analyzing sums involving Dirichlet characters. We then extend this method to characterize the distribution of powers in finite fields.
\end{abstract}
\indent Let $n>1$ be an integer. A \emph{Dirichlet character mod} $n$ is a homomorphism \newline$\chi:(\mathbb{Z}/n\mathbb{Z})^{\times} \to \mathbb{C}^{\times}$. In other words, $\chi$ is completely multiplicative and periodic mod $n$. It is conventional to treat $\chi$ as a function of integers and to set $\chi(a)=0$ if $\gcd(a,n)>1$ (preserving multiplicativity).

\indent From Euler's theorem, it follows that if $\gcd(a,n)=1$, then \begin{align*}\chi(a)^{\varphi(n)}&=\chi(a^{\varphi(n)})\\
&=\chi(1)\\
&=1,
\end{align*}
so the nonzero values of $\chi$ are all $\varphi(n)$th roots of unity. The \emph{principal character} is the character for which $\chi(a)$ is 1 if $\gcd(a,n)=1$ and 0 otherwise, which we write as $\chi_1$. Dirichlet characters have some nice orthogonality properties, including \[\sum_{\substack{1 \le a \le n\\ \gcd(a,n)=1}}\chi(a)=\begin{cases} \varphi(n)\;\text{if}\;\chi=\chi_1\\
\;\;\;0\;\;\;\text{otherwise},
\end{cases}\]
as well as \[\sum_{\chi}\chi(a)=\begin{cases} \varphi(n)\;\text{if}\;a\equiv 1 \pmod n\\
\;\;\;0\;\;\;\text{otherwise},
\end{cases}\]
where the sum is taken over all characters mod $n$. 
\newline
\indent Here, we will restrict ourselves to the case when $n$ is an odd prime. Doing so not only makes computations simpler but also allows us to use the fact that there exists a primitive root modulo any prime $p$, which will be denoted $g$. Since $g$ generates $(\mathbb{Z}/p\mathbb{Z})^{\times}$, it follows that any character mod $p$ is completely determined by its value at $g$. 
\newline
\indent In the work that follows, we consider $n$th power residues modulo a prime $p$. Note that if $\gcd(n,p-1)=d$, then there are integers $u,v$ such that $un+v(p-1)=d$, so $g^d=g^{un}$ is an $n$th power residue. Thus, the $d$th power residues are $n$th power residues, and vice versa. We therefore assume that $n \mid p-1$.
\newline
\indent Before proving any results, we need a key lemma.
\newline
\newline
\textbf{Lemma 1:} Let $\chi$ be a non-principal Dirichlet character mod $p$. Then \[\left|\sum_{n=0}^{p-1}\chi(n)e^{2\pi ni/p}\right|=\sqrt{p}.\]
\emph{Proof:} We have 
\begin{align*}\left|\sum_{n=0}^{p-1}\chi(n)e^{2\pi ni/p}\right|^2&=\left(\sum_{n=0}^{p-1}\chi(n)e^{2\pi ni/p}\right)\left(\sum_{n=0}^{p-1}\overline{\chi(m)}e^{-2\pi mi/p}\right)\\
&=\sum_{0 \le n,m\le p-1}\chi(n)\overline{\chi(m)}e^{2\pi (n-m)i/p}.
\end{align*}
Making the substitution $n=m+k$ gives
\[\sum_{0 \le n,m \le p-1}\chi(n)\overline{\chi(m)}e^{2\pi (m-n)i/p}=\sum_{0 \le k,m \le p-1} \chi(m+k)\overline{\chi(m)}e^{2\pi ki/p}.\]
Since $\chi(0)=0$, we can let $m$ be nonzero. This allows us to invert $m$, and so 
\[\sum_{0 \le k,m \le p-1} \chi(m+k)\overline{\chi(m)}e^{2\pi ki/p}=\sum_{0 \le k\le p-1}\sum_{1 \le m \le p-1} \chi(1+km^{-1})e^{2\pi ki/p}.\]
If $k=0$, then $\chi(1+km^{-1})=1$, so the inner sum is $p-1$. Otherwise, as $m$ varies, $1+km^{-1}$ varies over all elements of $\mathbb{Z}/p\mathbb{Z}$ except 1. The inner sum is therefore $-e^{2\pi ki/p}$ since $\chi\ne \chi_1$. It follows that 
\begin{align*}\left|\sum_{n=0}^{p-1}\chi(n)e^{2\pi ni/p}\right|^2&=p-1-\sum_{k=1}^{p-1}e^{2\pi ki/p}\\
&=p,
\end{align*}
and we are done. $\hfill\blacksquare$
\newline
\newline
\indent We now give our first theorem relating Fourier series to the distribution of $m$th power residues. For brevity, we let $R_m$ denote the set of $m$th power residues mod $p$ between 1 and $p-1$.
\newline
\newline
\textbf{Theorem 1:} Let $m>1$ be an integer and $p\equiv 1 \pmod m$ be a prime. Let $f:[0,1]\to \mathbb{R}$ be a function whose Fourier series converges pointwise to $f$ on $(0,1)$, say \[f(x)=\sum_{n=-\infty}^\infty a_n e^{2\pi nix}.\]
Suppose the sum \[S(f)=\sum_{n\ne 0}|a_n|\]
converges. Then \[\left|\sum_{k \in R_m} f\left(\frac{k}{p}\right)-\frac{1}{m}\sum_{k=1}^{p-1}f\left(\frac{k}{p}\right)\right|\le\left(1-\frac{1}{m}\right)S(f)\sqrt{p}.\]
\emph{Proof:} Let $\chi$ be a non-principal character mod $p$ for which $\chi(g)$ is an $m$th root of unity. Then \[\sum_{k=1}^{p-1}\chi(k)f\left(\frac{k}{p}\right)=\sum_{n=-\infty}^\infty a_n\left[\chi(1)e^{2\pi ni/p}+\chi(2)e^{4\pi ni/p}+\cdots+\chi(p-1)e^{(2p-2)\pi ni/p}\right].\]
Observe that the RHS can be rewritten as 
\[\sum_{n=-\infty}^\infty a_n\overline{\chi(n)}\left[\chi(n)e^{2\pi ni/p}+\chi(2n)e^{4\pi ni/p}+\cdots+\chi(n(p-1))e^{(2p-2)\pi ni/p}\right],\]
since $\chi(n)\overline{\chi(n)}=1$ unless $p \mid n$, in which case the summand is 0. This sum, in turn, is equal to 
\[\sum_{n=-\infty}^\infty a_n\overline{\chi(n)}\left[\chi(1)e^{2\pi i/p}+\chi(2)e^{4\pi i/p}+\cdots+\chi(p-1)e^{(2p-2)\pi i/p}\right].\]
The bracketed expression is precisely the Gauss sum from Lemma 1, so taking absolute values gives 
\begin{align*}\left|\sum_{k=1}^{p-1}\chi(k)f\left(\frac{k}{p}\right)\right|&=\sqrt{p}\left|\sum_{n=-\infty}^\infty a_n\overline{\chi(n)}\right|\\
& \le \sqrt{p}\sum_{n \ne 0}|a_n|\\
&=S(f)\sqrt{p}.
\end{align*}
We therefore have\[\left|\sum_{k=1}^{p-1}\chi(k)f\left(\frac{k}{p}\right)\right| \le S(f)\sqrt{p}.\]
Summing over all non-principal $\chi$ for which $\chi(g)^m=1$, we get 
\begin{align*} \sum_{\chi \ne 1}\left|\sum_{k=1}^{p-1}\chi(k)f\left(\frac{k}{p}\right)\right| \le (m-1)S(f)\sqrt{p}.
\end{align*}
Again by the triangle inequality, 
\begin{align*}\left|\sum_{\chi \neq 1}\sum_{k=1}^{p-1}\chi(k)f\left(\frac{k}{p}\right)\right|&\le (m-1)S(f)\sqrt{p}\\
\left|m\sum_{k \in R_m}f\left(\frac{k}{p}\right)-\sum_{k=1}^{p-1} f\left(\frac{k}{p}\right)\right|&\le(m-1)S(f)\sqrt{p},
\end{align*}
where in the second line we used the fact that summing over all $\chi$ for which $\chi(g)^m=1$ gives each $m$th power residue weight $m$ and each nonresidue weight 0. Dividing both sides of the inequality by $m$, we have the desired result. $\hfill\blacksquare$
\newline
\newline
\indent We will now use Theorem 2 to prove another result that says that the $m$th power residues are randomly distributed, roughly speaking. We require the following lemmas.
\newline
\newline
\textbf{Lemma 2:} Let $t,x$ be real numbers with $0<t\le 1$ and $0<x<1$. Then \[\sum_{n=1}^\infty \frac{t^n\sin 2\pi nx}{n}=\arctan\left(\frac{t \sin 2\pi x}{1-t\cos 2\pi x}\right).\]
\emph{Proof:} We have 
\begin{align*}\sum_{n=1}^\infty \frac{t^n\sin 2\pi nx}{n}&=\text{Im}\left(\sum_{n=1}^\infty \frac{t^ne^{2\pi nix}}{n}\right)\\
&=-\text{Im}[\log(1-te^{2\pi ix})]\\
&=\arctan\left(\frac{t\sin 2\pi x}{1-t\cos 2\pi x}\right).\end{align*}
The last step requires some care in choosing the branch of $\log z$, but it suffices to check that equality holds for a single pair $(t,x)$. $\hfill\blacksquare$
\newline
\newline
\textbf{Lemma 3:} Let $\delta,t$ be positive numbers with $\delta\le 1/3$ and $t \le 1/(1+3\delta^2)$. Then \[\frac{1-t\cos 2 \pi \delta}{t\sin 2\pi \delta}-\pi \delta \le \frac{\pi}{2\delta}\left(\frac{1}{t}-1\right).\]
\emph{Proof:} We have
\begin{align*}\frac{1-t\cos 2 \pi \delta}{t\sin 2\pi \delta}-\pi \delta&=\left(\frac{1}{t}-1\right)\csc 2\pi\delta+\csc2\pi\delta-\cot2\pi\delta-\pi\delta\\
&=\left(\frac{1}{t}-1\right)\csc 2\pi\delta+(\tan\pi\delta-\pi\delta).
\end{align*}
Now both $x \csc2\pi x$ and $(\tan\pi x)/x$ are increasing on $(0,1/3)$, so we have 
\[\csc 2\pi \delta \le \frac{2}{3\sqrt{3}\delta}, \tan\pi \delta \le 3\sqrt{3}\delta.\]
Thus, it suffices to prove that 
\[\frac{2}{3\sqrt{3}\delta}\left(\frac{1}{t}-1\right)+(3\sqrt{3}-\pi)\delta \le \frac{\pi}{2\delta}\left(\frac{1}{t}-1\right),\]
which after rearranging becomes 
\[\frac{1}{t}-1 \ge \frac{(3\sqrt{3}-\pi)\delta^2}{\frac{\pi}{2}-\frac{2}{3\sqrt{3}}}.\]
The RHS is less than $3\delta^2$, so taking $t \le 1/(1+3\delta^2)$ is sufficient and we are done.
\newline
\newline
\textbf{Theorem 2:} Let $m>1$ be an integer and $C$ be a constant greater than $\frac{3}{\pi}\left(1-\frac{1}{m}\right)$. Then for all sufficiently large primes $p \equiv 1 \pmod m$, the number of $m$th power residues in any interval $(a,b)\subset (0,p)$ is within $C\sqrt{p}\log p$ of $(b-a)/m$.
\newline
\newline
\emph{Proof:} Take a prime $p\equiv 1 \pmod m$ and let $\alpha=a/p,\beta=b/p$. Consider the function \[f(x)=\begin{cases} 1\;\text{if}\;\alpha<x<\beta\\
\frac{1}{2}\;\text{if}\;x=\alpha\;\text{or}\;x=\beta\\
0\;\text{otherwise}.
\end{cases}\]
It is clear that
\[\left|\sum_{k \in R_m}f\left(\frac{k}{p}\right)-\frac{1}{m}\sum_{k=1}^{p-1}f\left(\frac{k}{p}\right)\right|\]
is the difference between the number of $m$th power residues in $(a,b)$ and \newline $(b-a)/m$, up to some small constant. To bound this quantity, we examine the Fourier series of $f$. It is straightforward to find that \[f(x)=(\beta-\alpha)+\frac{1}{\pi}\sum_{n=1}^\infty \frac{\sin(2\pi n(x-\alpha)
)-\sin(2\pi n(x-\beta))}{n}.\]
This series converges pointwise to $f$, but \[\sum_{n \ne 0} |a_n|=\sum_{n \ne 0}\frac{|e^{-2\pi ni\alpha}-e^{-2\pi ni\beta}|}{2\pi |n|},\]
which may diverge. Thus, we cannot directly apply Theorem 2 and instead must approximate $f$; this will give us the error bound of $O(\sqrt{p} \log p)$.
\newline
\indent Consider the functions \[f_t(x)=(\beta-\alpha)+\frac{1}{\pi}\sum_{n=1}^\infty t^n \frac{\sin(2\pi n(x-\alpha)
)-\sin(2\pi n(x-\beta))}{n}\]
for $0<t<1$. All these functions satisfy the conditions of Theorem 2, so \[\left|\sum_{k \in R_m} f_t\left(\frac{k}{p}\right)-\frac{1}{m}\sum_{k=1}^{p-1}f_t\left(\frac{k}{p}\right)\right|\le \left(1-\frac{1}{m}\right)S(f_t)\sqrt{p}.\]
If we can make $|f(x)-f_t(x)|$ small by taking $t$ near 1, then we can obtain a similar bound for $f$. We therefore need to determine the rate at which $f_t$ converges to $f$.
\newline
\indent Let \[g_t(x)=\sum_{n=1}^\infty \frac{\sin 2\pi nx}{n}-\sum_{n=1}^\infty \frac{t^n \sin 2\pi nx}{n},\]
where $x$ is in some interval $[\delta,1-\delta]$ to avoid the discontinuities at $x=0$ and $x=1$. Differentiating with respect to $x$ yields 
\begin{align*}g'_t(x)&=-\pi-\frac{2\pi (t\cos 2\pi x-t^2)}{1-2t \cos 2\pi x+t^2}\\
&=\frac{-\pi(1-t^2)}{1-2t\cos 2\pi x+t^2}
\end{align*}
by Lemma 2. Since $t<1$, $g_t'$ is negative and so the extreme values of $g_t$ occur at the endpoints of the interval considered. Furthermore, $g_t(1-\delta)=-g_t(\delta)$, so 
\begin{align}|g_t(x)|\le |g_t(\delta)|&=\left|\frac{\pi-2\pi \delta}{2}-\arctan\left(\frac{t\sin 2\pi \delta}{1-t\cos 2\pi \delta}\right)\right|,
\end{align}
again by Lemma 2. From $0<\delta<1/2$ it follows that \[\frac{t \sin 2\pi \delta}{1-t\cos 2\pi \delta}>0.\] Thus, we can combine the $\pi/2$ and arctangent terms to get 
\[\left|\arctan\left(\frac{1-t\cos 2\pi \delta}{t\sin 2\pi \delta}\right)-\pi \delta\right|\]
for the bound on $|g_t|$. We remove the absolute value bars and use Laurent series to get
\begin{align*}|g_t(x)| &\le \arctan\left(\frac{1-t\cos 2\pi \delta}{t\sin 2\pi \delta}\right)-\pi \delta \\
&\le \frac{1-t\cos 2\pi \delta}{t\sin 2\pi \delta}-\pi \delta\\
&=\left(\frac{1}{t}-1\right)\cdot \frac{1}{2\pi \delta}+O(\delta).
\end{align*}
By increasing the constant $1/(2\pi)$ to $\pi/2$, we can ignore the higher-order terms for $\delta$ sufficiently small (depending on $t$). 
\newline
\indent By Lemma 3, if $\delta\le 1/3$ and $t\le 1/(1+3\delta^2)$, then
\begin{align*}\left|\sum_{n=1}^\infty \frac{\sin 2\pi nx}{n}-\sum_{n=1}^\infty \frac{t^n \sin 2\pi nx}{n}\right|\le \frac{\pi}{2\delta}\left(\frac{1}{t}-1\right)
\end{align*}
for $\delta \le x \le 1-\delta$. We can then bound $|f_t(x)-f(x)|$  as follows:
\begin{align*}|f(x)-f_t(x)|&=\frac{1}{\pi}\left|\sum_{n=1}^\infty \frac{(1-t^n)[\sin(2\pi n(x-\alpha))-\sin(2\pi n(x-\beta))]}{n}\right|\\
&\le \frac{1}{\pi}\left|\sum_{n=1}^\infty \frac{(1-t^n)\sin(2\pi n(x-\alpha))}{n}\right|+\frac{1}{\pi}\left|\sum_{n=1}^\infty \frac{(1-t^n)\sin(2\pi n(x-\beta))}{n}\right|\\
&\le \frac{1}{\delta}\left(\frac{1}{t}-1\right),
\end{align*}
where $\delta$ is chosen so that $x$ is at least a distance $\delta$ from the discontinuities at $\alpha$ and $\beta$. Since we care only about $x=1/p,2/p,\ldots,(p-1)/p$, the optimal $\delta$ is \[\delta=\min_{1 \le k \le p-1}\min\left\{\left|\frac{k}{p}-\alpha\right|,\left|\frac{k}{p}-\beta\right|\right\}.\]
We can do far better, however, by setting aside the two multiples of $1/p$ nearest to $\alpha$ and similarly for $\beta$ (if $\alpha$ or $\beta$ is a multiple of $1/p$, it does not matter which of the neighboring points we choose). Ignoring these four values of $k/p$ will increase our error bound by some small quantity. On the other hand, we can now safely take $\delta=1/p$ since the remaining values of $k/p$ are more than $1/p$ away from $\alpha$ and $\beta$. Then for $p$ sufficiently large, \[\left|f\left(\frac{k}{p}\right)-f_t\left(\frac{k}{p}\right)\right|\le p\left(\frac{1}{t}-1\right)\]
for all but the four special values of $k$, which we handle separately. From (1), it follows that $|g_t(x)| \le \pi/2$ for any $t,x$, so \[|f(x)-f_t(x)| \le \frac{1}{\pi}\left(\frac{\pi}{2}+\frac{\pi}{2}\right)=1.\]
Letting \[\epsilon=p\left(\frac{1}{t}-1\right),\]
we therefore have \[\left|f\left(\frac{k}{p}\right)-f_t\left(\frac{k}{p}\right)\right| \le \begin{cases}1\;\text{if}\;\frac{k}{p}\;\text{is one of four nearest to }\;\alpha,\beta\\
\epsilon\;\text{otherwise}.
\end{cases}\]
Then by Theorem 2 applied to $f_t$ and the triangle inequality, \[\left|\sum_{k \in R_m} f\left(\frac{k}{p}\right)-\frac{1}{m}\sum_{k=1}^{p-1}f\left(\frac{k}{p}\right)\right|\le\left(1-\frac{1}{m}\right)S(f_t)\sqrt{p}+\frac{2p-2}{m}\epsilon+4,\]
and we wish to choose $t$ so that the size of the the RHS is minimal. We take $\epsilon=1/\sqrt{p}$ so that the second and third terms are negligible compared to the first. This corresponds to \[t=\frac{1}{1+p^{-3/2}}.\]
The conditions of Lemma 3 are satisfied since $\delta=1/p\le 1/3$ and $t<1/(1+3\delta^2)$. Now we have 
\begin{align*}S(f_t)&=\sum_{n\ne 0}t^{|n|}\frac{|e^{-2\pi ni\alpha}-e^{-2\pi ni\beta}|}{2\pi |n|}\\
&\le \frac{2}{\pi}\sum_{n=1}^\infty \frac{t^n}{n}\\
&=-\frac{2}{\pi}\log(1-t)\\
&=\frac{2}{\pi}\log(p^{3/2}+1).
\end{align*}
Hence, the dominant term in the bound is at most \[\frac{2}{\pi}\left(1-\frac{1}{m}\right)\sqrt{p}\log(p^{3/2}+1).\]
Finally, we observe that the lower bound on $p$ needed to get these estimates depends on $C$, but not on $\alpha$ or $\beta$. The result immediately follows. $\hfill\blacksquare$
\newline
\newline
\indent Theorem 2 is, up to a constant, a consequence of the Polya-Vinogradov inequality [1]. Rather than looking at power residues in $\mathbb{Z}/p\mathbb{Z}$ over an interval, we will now examine the distribution of power residues in a field extension over a higher-dimensional box. In the work that follows, we write $\mathbb{F}_p$ in lieu of $\mathbb{Z}/p\mathbb{Z}$ to emphasize that we are working with field extensions. Here, $R_m$ will denote the set of $m$th powers in the chosen finite field.
\newline
\newline
\textbf{Theorem 3:} Fix integers $d,m>1$ and let $p$ be a prime for which $m \mid p^d-1$. Choose a polynomial of degree $d$ with integer coefficients that is irreducible over $\mathbb{F}_p$, and let $\xi$ be one of its roots. It follows that $\mathbb{F}_p(\xi)$ is a field, and every element can be written uniquely in the form \[c_0+c_1\xi+c_2\xi^2+\cdots+c_{d-1}\xi^{d-1},\]
where the $c_i$'s are in $\mathbb{F}_p$. Then for any $d$-dimensional box \[R=[a_0,b_0]\times [a_1,b_1]\times\cdots\times [a_{d-1},b_{d-1}]\subset(0,p)^d,\]
the number of $m$th powers in $\mathbb{F}_p(\xi)$ with $(c_0,c_1,\ldots,c_{d-1})\in R$ is, up to a small error, \[\frac{(\lfloor b_0\rfloor-\lceil a_0 \rceil)(\lfloor b_1\rfloor-\lceil a_1\rceil)\cdots(\lfloor b_{d-1}\rfloor-\lceil a_{d-1}\rceil)}{m}.\]
This error is bounded in absolute value by \[\left[\frac{2}{\pi}\sqrt{p}\log(3p^d+1)\right]^d+\frac{2d}{m}p^{d/2}.\]
\emph{Proof:} First note that since $\xi$ is a root of a polynomial that is irreducible over $\mathbb{F}_p$, the extension $\mathbb{F}_p(\xi)$ is indeed a field. Call this field $F$. The extension is of degree $d$, so each element can be written uniquely as a linear combination of $1,\xi,\xi^2,\ldots,\xi^{d-1}$. We now need to prove the claim about the distribution of the $m$th powers. 
\newline
\indent Since $F$ is a finite field, its group of units is cyclic; let $g$ be a generator of this group. Since $m \mid p^d-1$, the nonzero $m$th powers in the field are precisely the powers of $g^m$. We define a character $\chi$ to be a multiplicative function from $F$ to $\mathbb{C}^\times$, which is completely determined by its value at $g$. As before, we will count the $m$th powers via these characters. 
\newline
\indent Next, observe that we can move the vertices of $R$ slightly without changing the number of lattice points inside it. We can replace $a_i$ with $\lceil a_i\rceil -1/2$ and $b_i$ with $\lfloor b_i \rfloor +1/2$, and this does not alter the number of $m$th powers in the box or the main term in the estimate. Thus, without loss of generality, we suppose the $a_i$'s and $b_i$'s are half-integers. Let $\alpha_j=a_j/p$, $\beta_j=b_j/p$, and \[R'=[\alpha_0,\beta_0]\times\cdots\times [\alpha_{d-1},\beta_{d-1}].\]
Now that $R$ has been scaled down by a factor of $p$, we let $1_{R'}$ be the indicator function of $R'$. We have 
\[1_{R'}(x_0,x_1,\ldots,x_{d-1})=\prod_{j=0}^{d-1}1_{[\alpha_j,\beta_j]}(x_j),\]
and each term on the RHS has a Fourier series that converges pointwise to the function except at $\alpha_j$ and $\beta_j$. Let $f(x_0,x_1,\ldots,x_{d-1})$ be the Fourier series of $1_{R'}$ and $g_j$ be the Fourier series of $1_{[\alpha_j,\beta_j]}$. We also introduce families of functions $g_{j,t}$, whose Fourier coefficients are those of $g_j$ weighted by $t^{|n|}$. Similarly, $f_t$ is a weighted version of $f$, defined to be the product of the $g_{j,t}$'s. Thus, we can write 
\begin{align}&f_t(x_0,x_1,\ldots,x_{d-1}) \nonumber\\
&=\;\;\;\sum_{\mathclap{n_0,n_1,\ldots,n_{d-1}}}\;\;\;a_{0,n_0}a_{1,n_1}\cdots a_{d-1,n_{d-1}}\exp[2\pi i(n_0x_0+n_1x_1+\cdots+n_{d-1}x_{d-1})],
\end{align}
where $a_{j,k}$ is the coefficient of $e^{2\pi kix}$ in the Fourier series of $g_{j,t}$.  
\newline
\indent To count the $m$th powers, we first need to bound 
\[\sum_{0 \le c_0,c_1\ldots,c_{d-1}\le p-1}\chi(c_0+c_1\xi+\cdots+c_{d-1}\xi^{d-1})\zeta^{n_0c_0+n_1c_1+\cdots+n_{d-1}c_{d-1}},\]
where $\chi$ is a non-principal character, $\zeta$ is a primitive $p$th root of unity, and the $n_i$'s are arbitrary integers. Let \[\psi(c_0+c_1\xi+\cdots+c_{d-1}\xi^{d-1})=\zeta^{n_0c_0+n_1c_1+\cdots+n_{d-1}c_{d-1}},\]
which, not coincidentally, is a homomorphism from $(F,+)$ to $\mathbb{C}^\times$. We wish to prove that  \[\left|\sum_{z \in F} \chi(z)\psi(z)\right|=\begin{cases}0\;\text{if}\;p \mid n_0,n_1,\ldots,n_{d-1}\\
p^{d/2}\;\text{otherwise}.\end{cases}\]
The first case is easy to check, since then $\psi(z)$ is always 1. Otherwise, we write
\[\left|\sum_{z \in F} \chi(z)\psi(z)\right|^2=\sum_{z,w\in F} \chi(z)\psi(z)\overline{\chi(w)\psi(w)}.\] Now make the substitution $z=w+u$ to get 
\begin{align*}\sum_{z,w\in F} \chi(z)\psi(z)\overline{\chi(w)\psi(w)}&=\sum_{u,w\in F}\chi(w+u)\overline{\chi(w)}\psi(w+u)\overline{\psi(w)}\\
&=\sum_{u \in F}\sum_{w\in F^\times}\chi(1+uw^{-1})\psi(u).
\end{align*}
If $u=0$, then $1+uw^{-1}=1$ for all $w$. Otherwise, for fixed nonzero $u$, $1+uw^{-1}$ varies over all elements of the field except 1. Thus,
\begin{align*}\sum_{u\in F}\sum_{w\in F^\times}\chi(1+uw^{-1})\psi(u)&=(p^d-1)-\sum_{u\in F^\times}\psi(u)\\
&=p^d-\sum_{u\in F}\psi(u),
\end{align*}
and it is straightforward to check that the sum on the RHS is zero when $\psi$ is not identically 1. This proves the claim. 
\newline
\indent From this bound and equation (2), we immediately get that 
\begin{align*}\left|\sum_{0 \le c_0,c_1,\ldots,c_{d-1}\le p-1} \chi(c_0+c_1\xi+\cdots+c_{d-1}\xi^{d-1})f_t\left(\frac{c_0}{p},\frac{c_1}{p},\ldots,\frac{c_{d-1}}{p}\right)\right|&\\\le p^{d/2}\Bigg|\sum_{n_0,n_1,\ldots,n_{d-1}\ne 0} a_{0,n_0}a_{1,n_1}\cdots &a_{d-1,n_{d-1}}\Bigg|\\
=p^{d/2}\prod_{j=0}^{d-1}\sum_{n_j\ne 0} |a_{j,n_j}|,\;\;\;\;\;\;\;\;\;\;\;\;\;\;\;\;\;\;\;\;\;&
\end{align*}
where in the last line we split the sum into a product and used the triangle inequality. Summing over those non-principal $\chi$ for which $\chi(g)^m=1$ and dividing through by $m$, we get 
\begin{align*}\left|\sum_{z \in R_m} f_t\left(\frac{z}{p}\right)-\frac{1}{m}\sum_{z \in F} f_t\left(\frac{z}{p}\right)\right|\le \left(1-\frac{1}{m}\right)p^{d/2}\prod_{j=0}^{d-1}\sum_{n_j \ne 0} |a_{j,n_j}|,
\end{align*}
where $z/p$ denotes  the point $(c_0/p,c_1/p,\ldots,c_{d-1}/p)$. For brevity, we let \[S(g_{j,t})=\sum_{n_j \ne 0}|a_{j,n_j}|.\]
Then the above inequality becomes
\begin{align}\left|\sum_{z \in R_m} f_t\left(\frac{z}{p}\right)-\frac{1}{m}\sum_{z \in F} f_t\left(\frac{z}{p}\right)\right|\le \left(1-\frac{1}{m}\right)p^{d/2}\prod_{j=0}^{d-1}S(g_{j,t}).
\end{align}
\indent Suppose we choose $t$ so that $|g_{j,t}(x)-g_j(x)|\le\epsilon$ for some fixed $\epsilon>0$ and all $j=0,1,\ldots,d-1$. Using the triangle inequality and the fact that $|g_{j,t}(x)|\le 1$ for all $j,t,x$, we get 
\begin{align*}|f_t(x_0,x_1,\ldots,x_{d-1})-f(x_0,x_1,\ldots,x_{d-1})|&=\left|\prod_{j=0}^{d-1} g_{j,t}(x_j)-\prod_{j=0}^{d-1} g_j(x_j)\right|\\
&\le d\epsilon.
\end{align*}
Then when we approximate $f$ with $f_t$ on the $p^d-1$ points in question, we will get an error of order $p^d\epsilon$. We therefore take $\epsilon=p^{-d/2}$ so that this term is negligible compared to the other terms in the bound. We use Lemma 3 to find sufficient conditions for
\[|g_{j,t}(x)-g_j(x)|\le p^{-d/2}\]
to hold for $x$ in some interval. Taking \[\delta=\frac{p^{-d/2}}{3},\;t=\frac{1}{1+3\delta^2}=\frac{1}{1+p^{-d}/3},\]
the conditions of the lemma hold and so 
\[|g_{j,t}(x)-g_j(x)| \le \frac{1}{\delta}\left(\frac{1}{t}-1\right)=p^{-d/2},\] where $x$ is at least $\delta$ away from the discontinuities of $g_j$. 
\newline
\indent Since the vertices of $R$ were taken to have half-integer coordinates, every lattice point $(c_0,c_1,\ldots,c_{d-1})$ satisfies 
$|c_j-a_j|\ge 1/2$ and $|c_j-b_j| \ge 1/2$ for all $j$. But $\delta<1/(2p)$, so when we rescale by $1/p$, we get \[\left|\frac{c_j}{p}-\alpha_j\right|,\left|\frac{c_j}{p}-\beta_j\right|>\delta\]
for all $j$. In other words, $c_j/p$ is at least $\delta$ away from the discontinuities of $g_j$ and so 
\[\left|f_t\left(\frac{z}{p}\right)-f\left(\frac{z}{p}\right)\right|\le dp^{-d/2}\]
for all $z\in F$.
\newline
\indent We now combine all this information to get the desired result. For the above choice of $t$, we have 
\begin{align*} S(g_{j,t})&=\sum_{n \ne 0}t^{|n|}\left|\frac{e^{-2\pi ni\alpha_j}-e^{-2\pi ni\beta_j}}{2\pi ni}\right|\\
&\le \frac{1}{\pi}\sum_{n \ne 0}\frac{t^{|n|}}{|n|}\\
&=\frac{2}{\pi}\log\left(\frac{1}{1-t}\right)\\
&=\frac{2}{\pi}\log(3p^d+1).
\end{align*}
From (3), it then follows that 
\[\left|\sum_{z \in R_m} f_t\left(\frac{z}{p}\right)-\frac{1}{m}\sum_{z \in F} f_t\left(\frac{z}{p}\right)\right|\le \left(1-\frac{1}{m}\right)p^{d/2}\left(\frac{2}{\pi}\log(3p^d+1)\right)^d.\]
This, of course, can be weakened slightly by removing the factor of $1-1/m$, and then we get the first term in the error bound from the theorem. 
\newline
\indent We have chosen $t$ so that $|f_t-f| \le dp^{-d/2}$ for the points in question, so we have
\[\left|\sum_{z \in R_m} f\left(\frac{z}{p}\right)-\frac{1}{m}\sum_{z \in F} f\left(\frac{z}{p}\right)\right|\le \left(\frac{2}{\pi}\sqrt{p}\log(3p^d+1)\right)^d+\frac{2(p^d-1)}{m}\cdot dp^{-d/2}.\]
This can be weakened slightly to give the cleaner bound 
\[\left|\sum_{z \in R_m} f\left(\frac{z}{p}\right)-\frac{1}{m}\sum_{z \in F} f\left(\frac{z}{p}\right)\right| \le \left(\frac{2}{\pi}\sqrt{p}\log(3p^d+1)\right)^d+\frac{2d}{m}p^{d/2}.\]
This is precisely the result we sought, since $f=1_{R'}$ except on the boundary of $R'$, but $\partial R'$ does not contain any of the rescaled lattice points. This concludes the proof. $\hfill\blacksquare$

\begin{centering}\section*{References}
\end{centering}
\noindent [1] C. Pomerance, \emph{Remarks on the P\'olya-Vinogradov Inequality}, \emph{Integers} (Proceedings of the Integers Conference, October 2009), 11A (2011), Article 19, 11pp.

\end{document}